\documentclass[10pt,reqno]{amsart}
\usepackage{amssymb,amscd,amsbsy}
\usepackage{amsthm}
\newcommand{\lb}{\linebreak}

\renewcommand{\a}{\alpha}
\renewcommand{\b}{\beta}
\newcommand{\g}{\gamma}

\newcommand{\s}{\sigma}

\newcommand{\B}{{\mathcal B}}

\newcommand{\h}{{\mathcal H}}

\newcommand{\X}{{\mathcal X}}

\newcommand{\R}{{\Bbb R}}

\newcommand{\bs}{\boldsymbol}

\newcommand{\bS}{{\boldsymbol S}}

\newcommand{\rf}[1]{(\ref{#1})}

\newcommand{\df}{\stackrel{\mathrm{def}}{=}}

\newcommand{\trace}{\operatorname{trace}}

\newcommand{\const}{\operatorname{const}}

\newcommand{\eeq}{\end{equation}}
\newcommand{\beq}{\begin{equation}}
\newcommand{\bay}{\begin{eqnarray}}
\newcommand{\ba}{\begin{align*}}
\newcommand{\ea}{\end{align*}}
\newcommand{\ey}{\end{eqnarray}}
\newcommand{\bey}{\begin{eqnarray*}}
\newcommand{\eey}{\end{eqnarray*}}

\newcommand{\be}{\infty}

\newtheorem{thm}{\hspace{\parindent}Theorem}[section]

\theoremstyle{remark}

\newtheorem*{rem*}{Remark}

\newcommand\fM{\frak M}

\begin{document}

\newcommand{\vse}{\vspace{.2in}}

\title{Triple operator integrals in Schatten--von Neumann norms and functions of perturbed noncommuting operators}

\maketitle
\begin{center}
\Large
Aleksei Aleksandrov$^{\rm a}$, Fedor Nazarov$^{\rm b}$, Vladimir Peller$^{\rm c}$
\end{center}

\begin{center}
\footnotesize
{\it$^{\rm a}$St.Petersburg Branch, Steklov Institute of Mathematics, Fontanka 27, 191023 St-Petersburg, Russia\\
$^{\rm b}$Department of Mathematics, Kent State University, Kent, OH 44242, USA\\
$^{\rm c}$Department of Mathematics, Michigan State University, East Lansing, MI 48824, USA}
\end{center}

\newcommand{\mt}{{\mathcal T}}

\footnotesize

{\bf Abstract.} We study perturbations of functions $f(A,B)$ of noncommuting self-adjoint operators $A$ and $B$ that can be defined in terms of double operator integrals. We prove that if $f$ belongs to the Besov class $B_{\be,1}^1(\R^2)$, then we have the following Lipschitz type estimate in the Schatten--von Neumann norm $\bS_p$, $1\le p\le2$ norm:
$\|f(A_1,B_1)-f(A_2,B_2)\|_{\bS_p}\le\const(\|A_1-A_2\|_{\bS_p}+\|B_1-B_2\|_{\bS_p})$. However, the condition $f\in B_{\be,1}^1(\R^2)$ does not imply the Lipschitz type estimate in $\bS_p$ with $p>2$. The main tool is Schatten--von Neumann norm estimates for triple operator integrals.

\medskip

\begin{center}
{\bf\large Int\'egrales triples op\'eratorielles en normes de Schatten--von Nemann\\
et fonctions d'op\'erateurs perturb\'es ne commutant pas }
\end{center}

\medskip

{\bf R\'esum\'e.} Nous examinons les perturbations de fonctions $f(A,B)$ d'op\'erateurs auto-adjoints $A$ et $B$ qui ne commutent pas. Telles fonctions peuvent \^etre d\'efinies
en termes d'int\'egrales doubles op\'eratorielles. Pour $f$ dans l'espace de Besov
$B_{\be,1}^1(\R^2)$ nous obtenons l'estimation lipschitzienne en norme de Schatten--von Neumann $\bS_p$, $1\le p\le2$: $\|f(A_1,B_1)-f(A_2,B_2)\|_{\bS_p}\le\const(\|A_1-A_2\|_{\bS_p}+\|B_1-B_2\|_{\bS_p})$.
Par ailleurs, la condition $f\in B_{\be,1}^1(\R^2)$ n'implique pas l'estimation lipschitzienne en norme de $\bS_p$ pour $p>2$. L'outil principale est des estimations d'int\'egrales triples op\'eratorielles dans les normes de $\bS_p$.

\normalsize

\

\begin{center}
{\bf\large Version fran\c caise abr\'eg\'ee}
\end{center}

\medskip

Nous continuons d'examiner les propri\'et\'es de fonctions d'op\'erateurs 
auto-adjoints perturb\'es qui ne commutent pas. Dans \cite{ANP} nous \'etudions
des estimations du type lipschitzien pour les fonctions d'op\'erateurs auto-adjoints qui ne commutent pas. Si $A$ et $B$ son des op\'erateurs auto-adjouints qui ne commutent pas forc\'ement, on d\'efini la fonction $f(A,B)$ comme l'int\'egrale double op\'eratorielle 
$$
f(A,B)\df\iint f(x_1,x_2)dE_A(x_1)dE_B(x_2)
$$
si $f$ est un multiplicateur de Schur (voir \cite{BS}, \cite{Pe1} et \cite{AP}
 pour des informations sur les multiplicateurs de Schur et sur les int\'egrales doubles op\'eratorielles). Ici $E_A$ et $E_B$ sont les mesures spectrales de $A$ et $B$.

Nous avons d\'emontr\'e dans \cite{ANP} que si $f$ est une fonction de la classe
de Besov $B_{\be,1}^1(\R^2)$, $A_1$, $B_1$, $A_2$, $B_2$ sont des op\'erateurs auto-adjoints tels que $A_2-A_1\in\bS_1$ (classe trace) et
$B_2-B_1\in\bS_1$, alors
$$
\|f(A_1,B_1)-f(A_2,B_2)\|_{\bS_1}
\le\const\|f\|_{B_{\be,1}^1}\max\{\|A_1-A_2\|_{\bS_1},\|B_1-B_2\|_{\bS_1}\}.
$$
Par ailleurs, nous avons \'etabli dans \cite{ANP} que la condition 
$f\in B_{\be,1}^1(\R^2)$ n'implique pas l'estimation lipschitzienne en norme 
op\'eratorielle.

Dans cette note nous consid\'erons le m\^eme probl\`eme dans la norme de Schatten--von Neumann $\bS_p$. On se trouve que si $1\le p\le2$, l'in\'egalit\'e suivante est vrai:
$$
\|f(A_1,B_1)-f(A_2,B_2)\|_{\bS_p}
\le\const\|f\|_{B_{\be,1}^1}\max\{\|A_1-A_2\|_{\bS_p},\|B_1-B_2\|_{\bS_p}\}.
$$

Par ailleurs, nous \'etablissons dans cette note que si $p>2$, il n'y a pas de nombre positif $M$ pour lequel on ait
$$
\|f(A_1,B_1)-f(A_2,B_2)\|_{\bS_p}
\le M\|f\|_{B_{\be,1}^1}\max\{\|A_1-A_2\|_{\bS_p},\|B_1-B_2\|_{\bS_p}\}.
$$

Pour d\'emontrer les r\'esultats ci-dessus nous utilisons la formule
\begin{align}
\label{intform}
f(A_1,B_1)-f(A_2,B_2)=&
\iiint\frac{f(x_1,y)-f(x_2,y)}{x_1-x_2}
\,dE_{A_1}(x_1)(A_1-A_2)\,dE_{A_2}(x_2)\,dE_{B_1}(y)\nonumber\\[.2cm]
+&\iiint\frac{f(x,y_1)-f(x,y_2)}{y_1-y_2}
\,dE_{A_2}(x)\,dE_{B_1}(y_1)(B_1-B_2)\,dE_{B_2}(y_2)
\end{align}
qui \'etait \'etablie dans \cite{ANP} (voir la partie anglaise pour la d\'efinition d'int\'egrales triples op\'eratorielles) et nous obtenons les propri\'et\'es suivantes 
d'int\'egrales triples op\'eratorielles:

\medskip

{\bf Th\'eor\`eme.} {\it Supposons que $\Psi$ appartient au produit tensoriel de Haagerup 
$L^\be\!\otimes_{\rm h}\!L^\be\!\otimes_{\rm h}\!L^\be$
et
$$
W=\iiint\Psi(x_1,x_2,x_3)\,dE_1(x_1)T\,dE_2(x_2)R\,dE_3(x_3)
$$
Alors on a:
\newline
{\em(i)} si $p\ge2$, $T\in\B(\h)$ et $R\in\bS_p$, alors $W\in\bS_p$ et
$$
\|W\|_{\bS_p}\le\|\Psi\|_{L^\be\!\otimes_{\rm h}\!L^\be\!\otimes_{\rm h}\!L^\be}
\|T\|\cdot\|R\|_{\bS_p};
$$
\newline
{\em(ii)} si $p\ge2$, $T\in\bS_p$ et $R\in\B(\h)$, alors $W\in\bS_p$ et
$$
\|W\|_{\bS_p}\le\|\Psi\|_{L^\be\!\otimes_{\rm h}\!L^\be\!\otimes_{\rm h}\!L^\be}
\|T\|_{\bS_p}\|R\|;
$$
\newline
{\em(iii)} si $1/p+1/q\le1/2$, $T\in\bS_p$ et $R\in\bS_q$, alors $W\in\bS_r$
avec $1/r=1/p+1/q$ et
$$
\|W\|_{\bS_r}\le\|\Psi\|_{L^\be\!\otimes_{\rm h}\!L^\be\!\otimes_{\rm h}\!L^\be}
\|T\|_{\bS_p}\|R\|_{\bS_q}.
$$}

Par ailleurs, nous d\'emontrons que si $p<2$, alors (i) et (2) sont faux.

Le produit tensoriel de Haagerup est d\'efini dans la partie anglaise de cette note. Remarquons que les diff\'erences divis\'ees dans la formule \rf{intform} ne doivent pas appartenir au produit tensoriel de Haagerup $L^\be\!\otimes_{\rm h}\!L^\be\!\otimes_{\rm h}\!L^\be$ pour toutes le fonctions $f\in B_{\be,1}^1(\R^2)$. Cependant, elles appartiennent aux produits du type de Haagerup 
$L^\be\!\otimes_{\rm h}\!L^\be\!\otimes^{\rm h}\!L^\be$ et
$L^\be\!\otimes^{\rm h}\!L^\be\!\otimes_{\rm h}\!L^\be$ (voir la partie anglaise pour les d\'efinitions).

\medskip

\begin{center}
------------------------------
\end{center}

\setcounter{section}{0}
\section{\bf Introduction}

In this note we continue studying functions of noncommuting self-adjoint operators under perturbation. In \cite{ANP} we studied Lipschitz type estimates for functions of noncommuting pairs of self-adjoint operators. Recall that for (not necessarily commuting) self-adjoint operators $A$ and $B$, we considered in \cite{ANP} the functional calculus $f\mapsto f(A,B)$ defined as follows.  For the class of functions $f$ that are defined at least on the cartesian product $\s(A)\times\s(B)$ of the spectra of the operators and such that $f$ is a Schur multiplier with respect to the spectral measures $E_A$ and $E_B$ of $A$ and $B$ the operator $f(A,B)$ is defined by
$
f(A,B)\df\iint f(x,y)\,dE_A(x_1)\,dE_B(x_2).
$
We refer the reader to \cite{Pe1} and \cite {AP} for the definition of Schur multipliers and double operator integrals; note also that the theory of double operator integrals was developed by Birman and Solomyak \cite{BS}.

It was explained in \cite{ANP} that if $f$ is a function in the Besov space 
$B_{\be,1}^1(\R^2)$, then $f$ is a Schur multiplier with respect to $E_A$ and $E_B$
for arbitrary bounded self-adjoint operators $A$ and $B$ (we refer the reader to \cite{Pee} and \cite{APPS}) for an introduction to Besov spaces.

In \cite{ANP} we established the following Lipschitz type estimate in trace norm
for functions $f$ in $B_{\be,1}^1(\R^2)$:
$$
\|f(A_1,B_1)-f(A_2,B_2)\|_{\bS_1}
\le\const\|f\|_{B_{\be,1}^1}\max\{\|A_1-A_2\|_{\bS_1},\|B_1-B_2\|_{\bS_1}\}.
$$
On the other hand, it was shown in \cite{ANP} that there is no such Lipschitz type estimate for functions in $B_{\be,1}^1(\R^2)$ in the operator norm.

Note that earlier it was shown in \cite{Pe1} and \cite{Pe2} that functions $f$ on 
the real line $\R$ of class $B_{\be,1}^1(\R)$ are {\it operator Lipschitz}, i.e.,
$$
\|f(A)-f(B)\|\le\const\|A-B\|
$$
for arbitrary self-adjoint operators on Hilbert space and such Lipschitz type estimates also hold in the trace norm (as well as in all Schatten--von Neumann norms). Recall that not all Lipschitz functions are operator Lipschitz, this was first proved by Farforovskaya in \cite{F}. 

However, it turned out that the situation with H\"older functions is quite different. It was shown in \cite{AP} that if $f$ is a H\"older function of 
order $\a$, 
$0<\a<1$, then it is {\it operator H\"older of order} $\a$, i.e., 
$\|f(A)-f(B)\|\le\const\|A-B\|^\a$ for self-adjoint operators $A$ and $B$.

The results of \cite{Pe1}, \cite{Pe2}, and \cite{AP} were extended in \cite{APPS} to functions of normal operators (in other words, to functions of commuting pairs of self-adjoint operators) and in \cite{NP} to functions on $n$-tuples of commuting self-adjoint operators.

In this paper we study Lipschitz type estimates for functions of noncommuting self-adjoint operators in Schatten--von Neumann norms. We show that for functions in $B_{\be,1}^1(\R^2)$, Lipschitz type estimates hold in the Schatten--von Neumann norm of $\bS_p$ for $p\in[1,2]$. However, there are no Lipschitz type estimates for $p>2$. 

To obtain Lipschitz type estimates, we represent the difference
$f(A_1,B_1)-f(A_2,B_2)$ in terms of triple operator integrals. In \S~2 we study Schatten--von Neumann properties of triple operator integrals.
In \S~3 we state Lipschitz type estimates in the Sachatten--von Neumann norm
$\bS_p$ for $p\in[1,2]$, while in \S~4 we show that such Lipschitz type estimates do not hold for $p>2$. 

\

\section{\bf Triple operator integrals in Schatten--von Neumann norms}

\

Let $E_1$, $E_2$, and $E_3$ be spectral measures on measurable spaces 
$(\X_1,\fM_1)$, $(\X_2,\fM_2)$, and $(\X_3,\fM_3)$ on a Hilbert space $\h$. In \cite{Pe3} triple operator integrals
\bay
\label{Wtoi}
W=\int_{\X_1}\int_{\X_2}\int_{\X_3}\Psi(x_1,x_2,x_3)\,dE_1(x_1)T\,dE_2(x_2)R\,dE_3(x_3)
\ey
were defined for bounded linear operators $T$ and $R$ and for functions
$\Psi$ in the projective tensor product 
$L^\be(E_1)\hat\otimes L^\be(E_2)\hat\otimes L^\be(E_3)$. For such functions $\Psi$ 
the following holds:
\bay
\label{Impp}
T\in\B(\h),\quad R\in\bS_p.\quad p\ge1\quad\Longrightarrow\quad W\in\bS_p
\ey
\bay
\label{Imppq}
T\in\bS_p,\quad R\in\bS_q,\quad\frac1p+\frac1q\le1
\quad\Longrightarrow\quad W\in\bS_r,\quad\mbox{where}\quad
\frac1r=\frac1p+\frac1q.
\ey

Later in \cite{JTT} the definition of triple operator integrals of the form \rf{Wtoi} was extended to functions $\Psi$ in the {\it Haagerup tensor product}
$L^\be(E_1)\!\otimes_{\rm h}\!L^\be(E_2)\!\otimes_{\rm h}\!L^\be(E_3)$
of the spaces $L^\infty(E_j)$, $j=1,2,3$. It consists of functions $\Psi$
that admit a representation
\bay
\label{Htp}
\Psi(x_1,x_2,x_3)=\sum_{j,k\ge0}\a_j(x_1)\b_{jk}(x_2)\g_k(x_3),
\ey
where $\{\a_j\}_{j\ge0},~\{\g_k\}_{k\ge0}\in L^\be(\ell^2)$, and 
$\{\b_{jk}\}_{j,k\ge0}\in L^\be(\B)$. Here $\B$ is the space of infinite matrices that induce bounded linear operators on $\ell^2$; $\B$ is endowed with the operator norm. We refer the reader to \cite{Pi} for Haagerup tensor products. Moreover, the following inequality holds:
$$
\|W\|\le
\|\Psi\|_{L^\be\!\otimes_{\rm h}\!L^\be\!\otimes_{\rm h}\!L^\be}\|T\|\cdot\|R\|,
$$
where
$$
\|\Psi\|_{L^\be\otimes_{\rm h}\!L^\be\otimes_{\rm h}\!L^\be}\df\inf\left\{
\|\{\a_j\}_{j\ge0}\|_{L^\be(\ell^2)}\|\{\b_{jk}\}_{j,k\ge0}\|_{L^\be(\B)}
\|\{\g_k\}_{k\ge0}\|_{L^\be(\ell^2)}\right\},
$$
the infimum being taken over all representations of $\Psi$ in the form \rf{Htp}.

However, unlike in the case 
$\Psi\in L^\be(E_1)\hat\otimes L^\be(E_2)\hat\otimes L^\be(E_3)$, for functions 
$\Psi$ in the Haagerup tensor product 
$L^\be(E_1)\!\otimes_{\rm h}\!L^\be(E_2)\!\otimes_{\rm h}\!L^\be(E_3)$, the situation with implications \rf{Impp} and \rf{Imppq} is more complicated. We proved in \cite{ANP} that there exist a function $\Psi$ in
$L^\be(E_1)\!\otimes_{\rm h}\!L^\be(E_2)\!\otimes_{\rm h}\!L^\be(E_3)$,
a bounded linear operator $T$, and an operator $R$ of trace class such that
the triple operator integral \rf{Wtoi} does not belong to trace class $\bS_1$.

Nevertheless, it turns out that implications in \rf{Impp} and \rf{Imppq} hold 
under certain assumptions on $p$ and $q$
for an arbitrary function $\Psi$ in 
$L^\be(E_1)\!\otimes_{\rm h}\!L^\be(E_2)\!\otimes_{\rm h}\!L^\be(E_3)$.

\begin{thm}
\label{SNSp}
Let $\Psi\in L^\be(E_1)\!\otimes_{\rm h}\!L^\be(E_2)\!\otimes_{\rm h}\!L^\be(E_3)$.
Then the following holds:
\newline
{\em(i)} if $p\ge2$, $T\in\B(\h)$, and $R\in\bS_p$, then the triple operator integral in {\em\rf{Wtoi}} belongs to $\bS_p$ and
\bay
\label{boSp}
\|W\|_{\bS_p}\le\|\Psi\|_{L^\be\!\otimes_{\rm h}\!L^\be\!\otimes_{\rm h}\!L^\be}
\|T\|\cdot\|R\|_{\bS_p};
\ey
\newline
{\em(ii)} if $p\ge2$, $T\in\bS_p$, and $R\in\B(\h)$, then the triple operator integral in {\em\rf{Wtoi}} belongs to $\bS_p$ and
$$
\|W\|_{\bS_p}\le\|\Psi\|_{L^\be\!\otimes_{\rm h}\!L^\be\!\otimes_{\rm h}\!L^\be}
\|T\|_{\bS_p}\|R\|;
$$
\newline
{\em(iii)} if $1/p+1/q\le1/2$, $T\in\bS_p$, and $R\in\bS_q$, then
then the triple operator integral in {\em\rf{Wtoi}} belongs to $\bS_r$
with $1/r=1/p+1/q$ and
$$
\|W\|_{\bS_r}\le\|\Psi\|_{L^\be\!\otimes_{\rm h}\!L^\be\!\otimes_{\rm h}\!L^\be}
\|T\|_{\bS_p}\|R\|_{\bS_q}.
$$
\end{thm}

We will see in \S~4 that statements (i) and (ii) of Theorem \ref{SNSp} do not hold for $p\in[1,2)$.

To prove Theorem \ref{SNSp}, we first prove statemants (i) and (ii) and then use
complex interpolation of bilinear operators, see Theorem 4.4.1 in \cite{BL}.

In \cite{ANP} we established the following formula for $f(A_1,B_1)-f(A_2,B_2)$ in the case when \lb$f\in B_{\be,1}^1(\R^2)$ and the pair $(A_2,B_2)$ is a trace class perturbation of the pair $(A_1,B_1)$:
\begin{align}
\label{intf}
f(A_1,B_1)-f(A_2,B_2)=&
\iiint\frac{f(x_1,y)-f(x_2,y)}{x_1-x_2}
\,dE_{A_1}(x_1)(A_1-A_2)\,dE_{A_2}(x_2)\,dE_{B_1}(y)\nonumber\\[.2cm]
+&\iiint\frac{f(x,y_1)-f(x,y_2)}{y_1-y_2}
\,dE_{A_2}(x)\,dE_{B_1}(y_1)(B_1-B_2)\,dE_{B_2}(y_2).
\end{align}
However, the divided differences
$$
(x_1,x_2,y)\mapsto\frac{f(x_1,y)-f(x_2,y)}{x_1-x_2}\quad\mbox{and}\quad
(x,y_1,y_2)\mapsto\frac{f(x,y_1)-f(x,y_2)}{y_1-y_2}
$$
do not have to belong to the Haagerup tensor product 
$L^\be\!\otimes_{\rm h}\!L^\be\!\otimes_{\rm h}\!L^\be$ (this follows from Theorem 3.1 of \cite{ANP}). Nevertheless, we defined in \cite{ANP} Haagerup like tensor products
$L^\be\!\otimes_{\rm h}\!L^\be\!\otimes^{\rm h}\!L^\be$ and 
$L^\be\!\otimes^{\rm h}\!L^\be\!\otimes_{\rm h}\!L^\be$, 
defined triple operator integrals for such Haagerup like tensor products, and
proved that the first divided difference belongs to 
$L^\be\!\otimes_{\rm h}\!L^\be\!\otimes^{\rm h}\!L^\be$
while the second divided difference belongs to 
$L^\be\!\otimes^{\rm h}\!L^\be\!\otimes_{\rm h}\!L^\be$.

We are going to use the above integral representation in the case when the pair
$(A_2,B_2)$ is an $\bS_p$ perturbation of the pair $(A_1,B_1)$ for $p\in[1,2]$.

\medskip

{\bf Definition.}
{\it A function $\Psi$ is said to belong to the Haagerup-like tensor product 
$L^\be\!\otimes_{\rm h}\!L^\be\!\otimes^{\rm h}\!L^\be$ of the first kind if it admits a representation
\bay
\label{yaH}
\Psi(x_1,x_2,x_3)=\sum_{j,k\ge0}\a_j(x_1)\b_{k}(x_2)\g_{jk}(x_3)
\ey
with $\{\a_j\}_{j\ge0},~\{\b_k\}_{k\ge0}\in L^\be(\ell^2)$, and 
$\{\g_{jk}\}_{j,k\ge0}\in L^\be(\B)$. For a bounded linear operator $R$ and 
for an operator $T$ of class $\bS_p$, $1\le p\le2$, we define the triple operator integral
\bay
\label{toiHttp}
W=\iiint\Psi(x_1,x_2,x_3)\,dE_1(x_1)T\,dE_2(x_2)R\,dE_3(x_3)
\ey
as the following continuous linear functional on the Schatten--von Neumann class 
$\bS_{p'}$, $1/p'=1-1/p$, (on the class of compact operators if $p=1$):
\bay
\label{fko}
Q\mapsto
\trace\left(\left(
\iiint
\Psi(x_1,x_2,x_3)\,dE_2(x_2)R\,dE_3(x_3)Q\,dE_1(x_1)
\right)T\right).
\ey
}

The fact that the linear functional \rf{fko} is continuous is a consequence of inequality \rf{boSp}, which also implies the following estimate:
$$
\|W\|_{\bS_p}\le\|\Psi\|_{L^\be\otimes_{\rm h}\!L^\be\otimes^{\rm h}\!L^\be}
\|T\|_{\bS_p}\|R\|,
$$
where $\|\Psi\|_{L^\be\otimes_{\rm h}\!L^\be\otimes^{\rm h}\!L^\be}$ is the infimum of 
$$
\|\{\a_j\}_{j\ge0}\|_{L^\be(\ell^2)}\|\{\b_k\}_{k\ge0}\|_{L^\be(\ell^2)}
\|\{\g_{jk}\}_{j,k\ge0}\|_{L^\be(\B)}
$$
over all representations in \rf{yaH}.

Similarly, suppose that $\Psi$ belongs to the {\it Haagerup like tensor product
$L^\be\!\otimes^{\rm h}\!L^\be\!\otimes_{\rm h}\!L^\be$
of the second kind}, i.e., $\Psi$ admits a representation
$$
\Psi(x_1,x_2,x_3)=\sum_{j,k\ge0}\a_{jk}(x_1)\b_{j}(x_2)\g_k(x_3),
$$
where $\{\b_j\}_{j\ge0},~\{\g_k\}_{k\ge0}\in L^\be(\ell^2)$, 
$\{\a_{jk}\}_{j,k\ge0}\in L^\be(\B)$, $T$ is a bounded linear operator, and $R\in\bS_p$, $1\le p\le2$. Then
the continuous linear functional 
$$
Q\mapsto
\trace\left(\left(
\iiint
\Psi(x_1,x_2,x_3)\,dE_3(x_3)Q\,dE_1(x_1)T\,dE_2(x_2)
\right)R\right)
$$
on the class $\bS_{p'}$ determines an operator 
$$
W\df\iiint\Psi(x_1,x_2,x_3)\,dE_1(x_1)T\,dE_2(x_2)R\,dE_3(x_3)
$$
of class $\bS_p$. Moreover,
$$
\|W\|_{\bS_p}\le
\|\Psi\|_{L^\be\otimes^{\rm h}\!L^\be\otimes_{\rm h}\!L^\be}
\|T\|\cdot\|R\|_{\bS_p}.
$$

The following result can be deduced from Theorem \ref{SNSp}.

\begin{thm}
\label{ftHtp}
Let $\Psi\in L^\be\!\otimes_{\rm h}\!L^\be\!\otimes^{\rm h}\!L^\be$.
Suppose that $T\in\bS_p$ and $R\in\bS_q$, where
$1\le p\le2$, and $1/p+1/q\le1$. Then the operator $W$ in {\em\rf{toiHttp}} belongs to $\bS_r$, $1/r=1/p+1/q$, and
$$
\|W\|_{\bS_r}\le\|\Psi\|_{L^\be\otimes_{\rm h}\!L^\be\otimes^{\rm h}\!L^\be}
\|T\|_{\bS_p}\|R\|_{\bS_q}.
$$
\end{thm}

A similar result holds for triple operator integrals defined above for functions 
$\Psi$ in
$L^\be\!\otimes^{\rm h}\!L^\be\!\otimes_{\rm h}\!L^\be$.

\

\section{\bf Lipschitz type estimates in $\bs{\bS_p}$ with $\bs{p\le2}$}

\

Recall that in \cite{ANP} we established formula \rf{intf} for functions
$f\in B_{\be,1}^1(\R^2)$ and pairs of self-adjoint operators $(A_1,B_1)$ and
$(A_2,B_2)$ such that $(A_2,B_2)$ is a trace class perturbation of $(A_1,B_1)$.
Moreover, we proved in \cite{ANP} that the first divided difference in \rf{intf} belongs to 
$L^\be\!\otimes_{\rm h}\!L^\be\!\otimes^{\rm h}\!L^\be$
while the second divided difference belongs to 
$L^\be\!\otimes^{\rm h}\!L^\be\!\otimes_{\rm h}\!L^\be$.

The following theorem shows that the same is true if we replace trace norm with the norm in $\bS_p$ for $p\in[1,2]$. It can be deduced from Theorem \ref{SNSp} and formula \rf{intf}.

\begin{thm}
\label{LtSp}
Let $1\le p\le2$ and let $f\in B_{\be,1}^1(\R^2)$. Suppose that $(A_1,B_1)$ and
$(A_2,B_2)$ are pairs of self-adjoint operators such that $A_2-A_1\in\bS_p$
and $B_2-B_1\in\bS_p$.
Then
\bay
\label{liya}
\|f(A_1,B_1)-f(A_2,B_2)\|_{\bS_p}
\le\const\|f\|_{B_{\be,1}^1}(\|A_1-A_2\|_{\bS_p}+\|B_1-B_2\|_{\bS_p}).
\ey
\end{thm}

We have defined functions $f(A,B)$ for $f$ in 
$B_{\be,1}^1(\R^2)$ only for bounded self-adjoint operators $A$ and $B$. 
However, as in the case of trace class perturbations (see \cite{ANP}), formula \rf{intf} allows us to define the difference $f(A_1,B_1)-f(A_2,B_2)$ in the case when $f\in B_{\be,1}^1(\R^2)$ and the 
self-adjoint operators $A_1,\,A_2,\,B_1,\,B_2$ are possibly unbounded once we know that the pair $(A_2,B_2)$ is an $\bS_p$ perturbation of the pair $(A_1,B_1)$, 
$1\le p\le2$.
Moreover, inequality \rf{liya} also holds for such operators.

\

\section{\bf Lipschitz type estimates Lipschitz type estimates in $\bs{\bS_p}$ with $\bs{p>2}$}

\

Recall that we showed in \cite{ANP} that the condition $f\in B_{\be,1}^1(\R^2)$ does not imply Lipschitz type estimates in the operator norm for functions of pairs of (not necessarily commuting) self-adjoint operators.
It turns out that the same is true in the $\bS_p$-norms for $p>2$.

The main result of this section shows that unlike in the case of commuting pairs of self-adjoint operators, the condition $f\in B_{\be,1}^1(\R^2)$ does not imply 
Lipschitz type estimates in the norm of $\bS_p$ with $p>2$.

\begin{thm}
\label{p>2}
Let $p>2$.
There is no positive number $M$ such that 
$$
\|f(A_1,B_1)-f(A_2,B_2)\|_{\bS_p}
\le M\|f\|_{L^\be(\R^2)}(\|A_1-A_2\|_{\bS_p}+\|B_1-B_2\|_{\bS_p})
$$ 
for all bounded functions $f$ on $\R^2$ with Fourier transform supported in $[-1,1]^2$ and for all finite rank self-adjoint operators $A_1,\,A_2,\,B_1,\,B_2$.
\end{thm}

The proof of Theorem \ref{p>2} uses a modification of the construction given in \cite{ANP}.

We conclude the paper with a theorem that can be deduced from Theorem \ref{p>2}.

\begin{thm}
Let $1\le p<2$.
There are spectral measures $E_1$, $E_2$ and $E_3$ on Borel subsets of $\R$, a function $\Psi$ in the Haagerup tensor product 
$L^\be(E_1)\!\otimes_{\rm h}\!L^\be(E_2)\!\otimes_{\rm h}\!L^\be(E_3)$ and an operator $Q$ in $\bS_p$ such that
$$
\iiint\Psi(x_1,x_2,x_2)\,dE_1(x_1)\,dE_2(x_2)Q\,dE_3(x_3)\not\in\bS_p.
$$
\end{thm}

\medskip

\noindent
{\bf Acknowledgements}

\medskip

The research of the first author is partially supported by RFBR grant 14-01-00198, the research of the second author is partially supported by NSF grant DMS 126562, the research of the third author is partially supported by NSF grant DMS 1300924.

\end{document}